\documentclass[aps,prl,twocolumn,groupedaddress]{revtex4-2}

\usepackage{graphicx}
\usepackage{amssymb,amsthm,amsmath,bm,dsfont}
\usepackage{ifthen}
\usepackage{color}
\usepackage{relsize}
\usepackage[T1]{fontenc}
\usepackage{fancyhdr}
\usepackage{enumitem}
\usepackage{tocloft}
\usepackage[colorlinks=true]{hyperref}  
\hypersetup{linkcolor=blue}
\renewcommand{\vec}[1]{{\bm{#1}}}

\newcommand{\R}{\mathbb{R}}
\newcommand{\dee}{\mathrm{d}}
\newcommand{\eps}{\varepsilon}
\newcommand{\expect}[1]{\mathds{E}\!\left[#1\right]}
\newcommand{\var}[1]{\mathds{V}\!\left[#1\right]}
\newcommand{\F}[2]{\vec{F}_#1^#2(\bm{\xi}_0)}
\newcommand{\T}{\top}

\begin{document}


\title{Unifying Lyapunov exponents with probabilistic uncertainty quantification}


\author{Liam Blake, John Maclean and Sanjeeva Balasuriya}
\affiliation{School of Computer and Mathematical Sciences, University of Adelaide, Adelaide SA 5005, Australia}


\date{\today}

\begin{abstract}
The Lyapunov exponent is well-known in deterministic dynamical systems as a measure for quantifying chaos and detecting coherent regions in physically evolving systems.  In this Letter, we show how the Lyapunov exponent  can be unified with stochastic sensitivity (which quantifies the uncertainty of an evolving uncertain system whose initial condition is certain) within a finite time uncertainty quantification framework in which both the dynamics and the initial condition of a continuously evolving $ n $-dimensional state variable are uncertain.
\end{abstract}


\maketitle


\section{\label{sec:intro}Introduction}

Many physically evolving systems can be modelled as nonlinear dynamical systems for a high-dimensional state variable $ \vec{x}_t $ evolving in continuous time $ t $.  A simple example is the motion of fluid parcels, for which $ \vec{x}_t  \in \R^3 $ is a  position vector driven by
the Eulerian velocity field.  Systems modelled via spatially-discretized partial differential equations also fall within this framework, 
with each component of $ \vec{x}_t $ containing the value
of the variable (such as a temperature, chemical concentration or charge density) at a grid location.  Alternatively, in coupled component models, the state-vector $ \vec{x}_t $ may for example contain the pressure, volume and temperature of a thermodynamic system, 
or the $ \mathrm{CO}_2 $-concentration in each of atmosphere, land, ocean and biomass.  

Predicting $ \vec{x}_t $ into the future requires excellent knowledge
of the rules for its evolution, as well as precise knowledge of the current measurement
$ \vec{x}_0 $.  Both these processes inevitably possess {\em uncertainty}: respectively, {\em uncertainty in
the model dynamics} (phenomena not completely explained by the theory, numerical errors from data, resolution and discretization) and {\em uncertainty in the initial condition} $ \vec{x}_0 $ (resolution and knowledge of `nowcasting' data).  All modern models used in forecasting and prediction should ideally seek to incorporate these two types of uncertainty; an example in climate modelling is the idea of stochastic parameterization \cite{BernerEtAl_2017_StochasticParameterizationShort}.  

A recent article by Blake et~al \cite{BlakeMacleanBalasuriya_SDE2023}  is able to address the uncertainty quantification (UQ) of the state vector $ \vec{x}_t $ at any later finite time $ t > 0 $ due to both these types of uncertainties.  While 
the methods are related to  stochastic differential equations (SDEs) and their linearizations \cite[e.g.]{Sanz-AlonsoStuart_2017_GaussianApproximationsSmall,Blagoveshchenskii_1962_DiffusionProcessesDepending}, we show
in this Letter that this formulation \cite{BlakeMacleanBalasuriya_SDE2023} allows the amalgamation of two hitherto unconnected entities: stochastic sensitivity \cite{Balasuriya_2020_StochasticSensitivityComputable} and the Lyapunov exponent \cite{Abarbanel_RevModernPhys1993,ShaddenEtAl_2005_DefinitionPropertiesLagrangian}.  

`Stochastic sensitivity' ($ S^2 $) was introduced by Balasuriya \cite{Balasuriya_2020_StochasticSensitivityComputable} as the variance of the maximal projection of the deviation 
of a noisy trajectory from its deterministic counterpart at a later time, in the presence of small noise but with a 
deterministic initial condition.  This is a natural theoretical tool within the probabilistic UQ 
viewpoint \cite{Cheng_IEEEAutomaticaSinica_2023Short,MacDonaldCampbell_StatsSurveys_2021,BerryHarlim_SIAMJUncQuant_2015}, and has seen usage in identifying more robust regions within unsteady fluid flows \cite{Balasuriya_2020_StochasticSensitivityComputable,BadzaEtAl_2023_HowSensitiveAre}
An explicit
expression for $ S^2 $ was derived in two dimensions \cite{Balasuriya_2020_StochasticSensitivityComputable}, and bounds exist in higher dimensions \cite{KaszasHaller_2020_UniversalUpperEstimate}.  

Lyapunov exponents, on the other hand, do not immediately appear to be connected to UQ. This purely deterministic concept quantifies the maximal stretching of an infinitesimal sphere centered at a fixed initial condition.  
Highly-used in detecting chaos \cite[viz.]{Abarbanel_RevModernPhys1993} because they quantify `sensitivity to initial conditions,' there
is nonetheless no formal connection to UQ in the sense of probability distributions.  Lyapunov exponents computed over a finite time, when considered as a field over all initial conditions, 
have seen tremendous use in detecting coherent regions and separators between them in fluids/geophysics \cite{ShaddenEtAl_2005_DefinitionPropertiesLagrangian,BalasuriyaKalampattelOuellette_Hyperbolic_2016}.  Suitable modifications of Lyapunov exponents continue to be used extensively across many other areas of physics: e.g., quantum \cite{BergamascoPhysRevE2023,XuScaffidiCao_PhysRevLett_2020,MaldacenaSandfordPHysRevD2016,MeldacenaShenkerStanford_JHighEnergyPhys_2016}, 
cosmology \cite{Staelens_PhysRevD_2023,CardosoPhysRevD2009}, electromagnetic \cite{Yamaguchi_PhysRevB_2023,JalabertPastawskiPhysRevLett2001,Borgogno_PhysPlasma_2011}, and 
statistical/thermal \cite{Tsallis_PhilTransRSoc_2023,Zou_AIP_2023,MeldacenaShenkerStanford_JHighEnergyPhys_2016}.
Currently, connections of Lyapunov exponents to stochastics is heuristic: comparison with correlations
in quantum chaos \cite{BergamascoPhysRevE2023,XuScaffidiCao_PhysRevLett_2020,Rozenbaum_PhysRevLett2017,MaldacenaSandfordPHysRevD2016,MeldacenaShenkerStanford_JHighEnergyPhys_2016} (based on semi-classical expectations, but displaying differences), and numerical estimation of stochastic flow barriers \cite{BalasuriyaGottwald_2018_EstimatingStableUnstable,GuoEtAl_2016_FiniteTimeLyapunovExponents}. 

In this Letter, we develop a UQ-formalism which  positions Lyapunov exponents within the probability framework.  We first describe modifications to Blake et~al \cite{BlakeMacleanBalasuriya_SDE2023} in the next Section.  We show in the next two Sections how these specialize to $ S^2 $  and to finite-time Lyapunov exponents respectively, thereby
unifying these concepts.

\section{\label{sec:uncertainty}{Relevant uncertainty results}}

This Section contains modifications and slight extensions to Blake et~al \cite{BlakeMacleanBalasuriya_SDE2023}.  The evolution of
a state vector $ \vec{x}_t \in \R^n $ from an initial time $ 0 $ to a general time $ t \in [0,T] $ with $ T < \infty $ is assumed to be governed by
\begin{equation}
\dee \vec{x}_t = \vec{u}\left( \vec{x}_t,t \right) \dee t + \eps \vec{\sigma} \left( \vec{x}_t, t \right) \dee \vec{W}_t \, ; \; 
\vec{x}_0 = \vec{\xi} \, , 
\label{eq:uncertain}
\end{equation}
where $ \vec{u} $ is the best-available deterministic dynamics (possibly given by data on a spatiotemporal grid),  and the It\^{o} form \cite{ito_NagoyaMathJ1950} is assumed for the SDE.
The {\em uncertainty in the dynamics} is captured via 
the presence of $ \vec{W}_t $, the classical $ n $-dimensional Wiener process, modulated by the state- and time-dependent
$ n \times n $ matrix $ \vec{\sigma} $.  The model noise is therefore multiplicative, and scaled by $
\eps $ where $ 0 < \eps \ll 1 $.  The {\em initial-condition uncertainty} is reflected through the fact that $ \vec{\xi} $ is sourced from the
Gaussian distribution $ \mathcal{N} \left( \vec{\xi}_0, \delta^2 \vec{\Xi}_0 \right) $, where the mean $ \vec{\xi}_0 \in \R^n $
is constant, and the matrix $ \vec{\Xi}_0 \in \R^{n \times n} $ is constant, positive-definite and symmetric, such that 
covariance $ \delta^2  \vec{\Xi}_0 $ has the scaling $ \delta $ where $ 0 < \delta \ll 1 $.  Thus, $ \eps $ and $ \delta $
respectively provide measures of the model and initial-condition uncertainties, each following a natural Gaussian interpretation.
When the expectation $ \expect{\centerdot} $ and variance $ \var{\centerdot} $ terminology is
used, these are with respect to the joint event space of the independent random variables
$ \vec{W}_t $ and $ \vec{\xi} $.

An equivalent understanding of the solution of \eqref{eq:uncertain} is that the distribution $ \rho(\vec{x},t) $ of $ \vec{x}_t $ obeys
the Fokker--Planck equation  
\begin{equation}
\frac{\partial \rho}{\partial t} \! + \! \vec{\nabla} \cdot \left( \rho \bm{u} \right) = \! \frac{\eps^2}{2} \vec{\nabla} \cdot \vec{\nabla} \cdot \left( \rho \, \vec{\sigma} \vec{\sigma}^\top \right) \, , 
\label{eq:fokkerplanck}
\end{equation}
with initial condition $ \rho(\vec{x},0) $ being the probability density function of $ {\mathcal N} \left( \vec{\xi}_0, \delta^2 \vec{\Xi}_0 \right) $.

If {\em both} uncertainties are turned off ($ \eps = 0 = \delta $), the fully uncertain model \eqref{eq:uncertain} becomes purely
deterministic:
\begin{equation}
\dee \vec{c}_t = \vec{u}\left( \vec{c}_t,t \right) \dee t  \, ; \; 
\vec{c}_0 = \vec{\xi}_0 \, . 
\label{eq:certain}
\end{equation}
Solutions to this {\em certain} model  will be denoted with the `flow-map' notation $ \vec{c}_t = \F{0}{t} $, the location to which
$ \vec{\xi}_0 $ is pushed from time $ 0 $ to time $ t $ 
by the deterministic flow \eqref{eq:certain}.  

In seeking validation of a commonly used practice \cite{Sanz-AlonsoStuart_2017_GaussianApproximationsSmall,Blagoveshchenskii_1962_DiffusionProcessesDepending,LawEtAl_2015_DataAssimilationMathematical}, Blake et~al \cite{BlakeMacleanBalasuriya_SDE2023} study the formal linearization of \eqref{eq:uncertain} with respect to solutions perturbed around
the deterministic solution.  The linearized (but still random) solution $ \vec{l}_t $ then satisfies
\begin{widetext}
\begin{equation}
\dee \vec{l}_t = 
\left[ \vec{u}\left( \F{0}{t},t \right) + \vec{\nabla} \vec{u} \left( \F{0}{t},t \right)  \left(
\vec{l}_t - \F{0}{t}\right) \right] \dee t + \eps \vec{\sigma} \left( \F{0}{t}, t \right) \dee \vec{W}_t \, ; \; 
\vec{l}_0 = \vec{\xi} \, , 
\label{eq:linear}
\end{equation}
\end{widetext}
subject to the same realization of the randomness
$ (\vec{W}_t, \vec{\xi}) $ as used in \eqref{eq:uncertain}.   Both $ \vec{l}_t $ and
$ \vec{x}_t $
depend on the two uncertainty parameters, and when important we write $ \vec{l}_t(\eps,\delta) $ and $ \vec{x}_t(\eps,\delta) $.
Now, 
the assumptions of Blake et~al \cite{BlakeMacleanBalasuriya_SDE2023} can be stated as
\begin{eqnarray*}
& & \left\| \vec{u} \left( \vec{x}, t \right) \right\| + \left\| \vec{\sigma} \left( \vec{x}, t \right) \right\|  \le U_0 \left( 1 + \left\| \vec{x} \right\| \right) \, , \\ 
& & \left\| \vec{\nabla} \vec{u} \left( \vec{x}, t \right) \right\| \le U_1 \, , \,
\left\| \vec{\nabla} \vec{\nabla} \vec{u} \left( \vec{x}, t \right) \right\| \le U_2 \, , \\
& & \left\| \vec{\sigma} \left( \vec{x}, t \right) \right\| \le W_0 \, \, \, \mathrm{and} \, \, \, 
\left\| \vec{\nabla} \vec{\sigma} \left( \vec{x}, t \right) \right\| \le W_1 
\end{eqnarray*}  
independent of $ (\vec{x}, t) \in \R^n \times [0,T] $.  The standard Euclidean norm is used here for vectors, and for higher-ranked tensors, 
the norm when used in this Letter is the operator norm induced by the tensor of lower rank.
The main theorem of Blake et~al \cite{BlakeMacleanBalasuriya_SDE2023} can, in the language of the current paper, be modified to: for any $ r \ge 1 $, 
\begin{equation}
\expect{\left\| \vec{x}_t  - \vec{l}_t \right\|^r}  \le  A \, \eps^{2r} + B \, \delta^r \, \eps^r
+  C \, \delta^{2r} \, ,
\label{eq:bound}
\end{equation}
in which $ A $, $ B $ and $ C $ do not depend on $ \eps $, $ \delta $, $ \vec{\xi}_0 $ or $ t $,
 but depend on $ U_{0,1,2} $, $ W_{0,1} $, $ n $ and $ \vec{\Xi}_0 $.
 Obtaining \eqref{eq:bound}
from the result of Blake et~al \cite{BlakeMacleanBalasuriya_SDE2023} requires some manipulations, including establishing that
the $ r $th moment of $ \vec{\xi} - \vec{\xi}_0 $ in relation to the distribution $ \mathcal{N} \left( \vec{\xi}_0,\delta^2 \vec{\Xi}_0 \right) $ is proportional to $ \delta^r  $, as can
be seen by working with the probability density function 
\footnote{The result arises by using the whitening transformation $ \vec{\nu} = \vec{\Xi}_0^{-1} \left( \vec{\xi} - \vec{\xi}_0 \right) / \sqrt{2 \delta^2} $ on the Gaussian variable $ \vec{\xi} $.}.
 Hence, we have the strong statement \eqref{eq:bound} on
how the true random solution $ \vec{x}_t $ to \eqref{eq:uncertain} approaches the solution $ \vec{l}_t $ to the
linearized equation \eqref{eq:linear} in the sense of expectation of all moments, as the noise parameters $ \left( \eps, \delta \right) \rightarrow \vec{0} $.  
Blake et~al \cite{BlakeMacleanBalasuriya_SDE2023} also derive an exact result: the strong solution to the linearized equation \eqref{eq:linear} is Gaussian with
\begin{equation}
\vec{l}_t \sim {\mathcal N} \left( \F{0}{t}, \vec{\Lambda}_t(\vec{\xi}_0; \eps, \delta) \right)\, , 
\label{eq:solution}
\end{equation}
where the covariance matrix is 
\begin{widetext}
\begin{equation}
\vec{\Lambda}_t (\vec{\xi}_0; \eps, \delta) :=  \delta^2  \, \vec{\nabla} \F{0}{t} \, \vec{\Xi}_0 \left[ \vec{\nabla} \F{0}{t} \right]^\T 
+ \eps^2 \vec{\nabla} \F{0}{t} \left[ 
\int_0^t \! \! M \left(\vec{\xi}_0,\tau \right) M \left(\vec{\xi}_0,\tau \right)^\T  \! \! \dee \tau \right] 
\left[ \vec{\nabla} \F{0}{t} \right]^\T
\label{eq:covariance}
\end{equation}
\end{widetext}
in which
\[
M(\vec{\xi}_0,\tau) := \left[ \vec{\nabla} \F{0}{\tau} \right]^{-1} \vec{\sigma} \left( \F{0}{\tau}, \tau \right) \, .
\]
By \eqref{eq:bound}, the distribution $ \rho $ for $ \vec{x}_t $ in \eqref{eq:fokkerplanck} is
well-approximated by the Gaussian probability density function \eqref{eq:solution} for small $ \eps $ and $ \delta $.

\section{\label{sec:s2}{Certain initial conditions and $ S^2 $}}

We quickly consider the case where the initial condition is assumed certain, which is the setting of stochastic sensitivity $ S^2 $, originally defined by Balasuriya \cite{Balasuriya_2020_StochasticSensitivityComputable} as the supremum over all projections of the variance of the scaled deviation from the deterministic solution, i.e., 
\begin{equation}
S^2(\vec{\xi}_0,t) := \lim_{\eps \downarrow 0} \! \sup_{\stackrel{\vec{p}}{\left\| \vec{p} \right\| = 1}} \! \var{ \vec{p}^{\T} \!\left( \frac{
\vec{x}_t(\eps,0) \! - \! \F{0}{t}}{\eps} \right) } \, . 
\label{eq:s2def}
\end{equation}
Since $ \vec{\xi} = \vec{\xi}_0 $ is considered {\em certain}, $ \delta = 0 $.
The original paper \cite{Balasuriya_2020_StochasticSensitivityComputable} was able to derive an analytical expression based on $ \vec{u} $ and $ \vec{\sigma} $ only
in two-dimensions.  While a {\em bound} was found in $ \R^n $ \cite{KaszasHaller_2020_UniversalUpperEstimate},  Blake et~al \cite{BlakeMacleanBalasuriya_SDE2023} showed that an 
analytic expression was possible.  Setting $ \delta = 0 $
in \eqref{eq:bound} and establishing that the limit $ \eps \rightarrow 0 $ of
$ \var{\vec{x}_t} / \eps^2 $ and $ \var{\vec{l}_t}/\eps^2 $ are identical, $ S^2 $ emerges as a special case of the
covariance \eqref{eq:covariance} as \cite{BlakeMacleanBalasuriya_SDE2023} 
\begin{equation}
S^2(\vec{\xi}_0,  t) = \frac{1}{\eps^2} \left\| \vec{\Lambda}_t \left( \vec{\xi}_0; \eps, 0 \right)  \right\|  \, .
\label{eq:s2}
\end{equation}
 In the above and subsequently, we remark that the operator norm $ \left\| \vec{A} \right\| $ of a square matrix $ \vec{A} $ is 
 easily computed as the square-root of the largest eigenvalue of $ \vec{A}^\top \vec{A} $.

\section{\label{sec:ftle}{Certain dynamics and Lyapunov exponents}}

While $ S^2 $'s connection to the covariance \eqref{eq:covariance} is anticipated given its stochastic
interpretation, the Lyapunov exponent is based on a {\em deterministic} model.  As our main result, we now establish that this too is a special case of the covariance \eqref{eq:covariance} for finite times.
For the purely deterministic system \eqref{eq:certain} the `finite-time Lyapunov exponent' (FTLE) is \cite{ShaddenEtAl_2005_DefinitionPropertiesLagrangian}
\begin{eqnarray}
\tilde{\lambda}(\vec{\xi}_0,t) & := & \frac{1}{t} \ln \sup_{\vec{\delta \xi}} \lim_{\left\| \vec{\delta \xi} \right\| \rightarrow \vec{0}} 
\frac{ \left\| \vec{F}_0^t \left( \vec{\xi}_0 + \vec{\delta \xi} \right) - \F{0}{t} \right\|}{ \left\| \vec{\delta \xi} \right\|} \nonumber \\
& = & \frac{1}{t} \ln \sup_{\vec{\delta \xi}} \frac{\left\| \vec{\nabla} \F{0}{t} \vec{\delta \xi} \right\|}{\left\| \vec{\delta \xi} \right\|}  
\nonumber \\
& = & \frac{1}{t} \ln \left\| \vec{\nabla} \F{0}{t} \right\| \, .
\label{eq:ftledef}
\end{eqnarray}
If an infinitesimal ball with radius vector $ \vec{\delta \xi} $ positioned at $ \vec{\xi}_0 $ at time $ 0 $, the flow from time $ 0 $ to $ t $ pushes it
to become an infinitesimal ellipse, and its `maximal stretching' $ \left\| \vec{\nabla} \F{0}{t} \right\| $  is the ratio of the longest axis of the ellipse to
the original radius.  The FTLE converts this to an exponential-in-time stretching rate, and is often viewed as a field
over $ \vec{\xi}_0 $ to determine initial conditions which are prone to most stretching or chaos in physical systems \cite{Abarbanel_RevModernPhys1993,BalasuriyaKalampattelOuellette_Hyperbolic_2016}, as well as a tool for determining
flow separators in Lagrangian coherent structures \cite{BalasuriyaEtAl_2018_GeneralizedLagrangianCoherent,HadjighasemEtAl_2017_CriticalComparisonLagrangian,ShaddenEtAl_2005_DefinitionPropertiesLagrangian}.   

Consider now the general stochastic model \eqref{eq:uncertain}. Comparing with \eqref{eq:s2def}, we define the uncertainty measure
\begin{equation}
Q^2(\vec{\xi}_0,t) := \lim_{\delta \downarrow 0} \! \sup_{\stackrel{\vec{p}}{ \left\| \vec{p} \right\| = 1}} \! \var{ \vec{p}^{\T} \! \left( \frac{
\vec{x}_t(0,\delta) \! - \! \F{0}{t}}{\delta} \right) } \, 
\label{eq:q2def}
\end{equation}
by inverting the roles of $ \eps $ and $ \delta $.
To proceed, we note \footnote{The trace of a positive definite matrix is the sum of its eigenvalues (all of which are positive), and
the norm of such a matrix is its largest eigenvalue.} that for any $ \vec{w} \in \R^n $,
\begin{equation}
\left\| \vec{w} \vec{w}^\T \right\| \le \mathrm{Trace} \left( \vec{w} \vec{w}^\T \right) = \left\| \vec{w} \right\|^2 \, .
\label{eq:trace}
\end{equation}
Now, since \eqref{eq:solution} with $ \eps = 0 $ tells us that
\[
\vec{l}_t(0,\delta) \sim {\mathcal N} \left( \F{0}{t} , \delta^2  \, \vec{\nabla} \F{0}{t} \, \vec{\Xi}_0 \left[ \vec{\nabla} \F{0}{t} \right]^\T 
\right) \, ,
\]
we have $ \expect{ \left(  \vec{l}_t(0,\delta) -  \F{0}{t} \right) / \delta } = \vec{0} $ and 
\[
\var{ \frac{ \vec{l}_t(0,\delta) -  \F{0}{t}}{\delta} } = 
\vec{\nabla} \F{0}{t} \, \vec{\Xi}_0 \left[ \vec{\nabla} \F{0}{t} \right]^\T \, .
\]
Upon defining
\[
\vec{w} := \frac{ \vec{x}_t (0,\delta)- \vec{l}_t(0,\delta)}{\delta} \, ,
\]
and applying $ \eps = 0 $ in the bound \eqref{eq:bound} with the choices $ r = 1 $ and $ r = 2 $ yields
\begin{equation}
 \expect{ \left\| \vec{w}  \right\|} \le C \delta \quad \mathrm{and} \quad
\expect{ \left\| \vec{w}  \right\|^2} \le C \, \delta^2 \, .
\label{eq:wdelta}
\end{equation}
Hence the norm of the covariance matrix of $ \vec{w} $ satisfies
\begin{eqnarray*}
\lim_{\delta \downarrow 0} \left\| \var{ \vec{w} } \right\| & =  & \lim_{\delta \downarrow 0} \left\|\expect{  \vec{w} \vec{w}^\T} - 
 \expect{ \vec{w}} \expect{ \vec{w}^\T}  \right\| \\ 
& \le & \lim_{\delta \downarrow 0} \left\|\expect{  \vec{w} \vec{w}^\T  } \right\| + \lim_{\delta \downarrow 0}  \left\| \expect{\vec{w}}
\lim_{\delta \downarrow 0} \expect{ \vec{w}^\T} \right\| \\
& \le & \lim_{\delta \downarrow 0} \expect{ \left\| \vec{w} \right\|^2} + 0  \\
& = & 0 \, , 
\end{eqnarray*}
by virtue of \eqref{eq:wdelta} and \eqref{eq:trace}.  Since $ \vec{w} $ then has expectation and variance $ 0 $ in the 
limit  $ \delta \downarrow 0 $, we are now in a position to compute $ Q^2 $ in \eqref{eq:q2def} as
\begin{eqnarray*}
Q^2(\vec{\xi}_0,t) & = & \lim_{\delta \downarrow 0} \! \sup_{\stackrel{\vec{p}}{ \left\| \vec{p} \right\| = 1}} \! \var{ \vec{p}^{\T} \! \left( \frac{
\vec{l}_t(0,\delta) \! - \! \F{0}{t}}{\delta} + \vec{w} \right) } \\
& = & \lim_{\delta \downarrow 0} \! \sup_{\stackrel{\vec{p}}{ \left\| \vec{p} \right\| = 1}}  \vec{p}^{\T} \! \var{ \! \left( \frac{
\vec{l}_t(0,\delta) \! - \! \F{0}{t}}{\delta} \right) } \vec{p} \\
& = & \! \sup_{\stackrel{\vec{p}}{ \left\| \vec{p} \right\| = 1}}  \vec{p}^{\T} 
\vec{\nabla} \F{0}{t} \, \vec{\Xi}_0 \left[ \vec{\nabla} \F{0}{t} \right]^\T \! \! \vec{p} \\
& = & \sup_{\stackrel{\vec{p}}{ \left\| \vec{p} \right\| = 1}}  \vec{p}^\T \vec{\nabla} \F{0}{t} \, \vec{\Psi}_0 \vec{\Psi}_0^\T \left[ \vec{\nabla} \F{0}{t} \right]^\T \! \! \vec{p} \\
& = & \sup_{\stackrel{\vec{p}}{ \left\| \vec{p} \right\| = 1}} \left\| \vec{\Psi}_0^\T \left[ \vec{\nabla} \F{0}{t} \right]^\T \! \! \vec{p} \right\|^2 \\
& = &  \left\| \vec{\Psi}_0^\T \left[ \vec{\nabla} \F{0}{t} \right]^\T \right\|^2 
=   \left\|  \vec{\nabla} \F{0}{t} \vec{\Psi}_0  \right\|^2 \, , 
\end{eqnarray*}
where we have used the Cholesky decomposition
\begin{equation}
\vec{\Xi}_0 = \vec{\Psi}_0 \vec{\Psi}_0^\T
\label{eq:cholesky}
\end{equation}
for positive definite symmetric $ \vec{\Xi}_0 $, in which $ \vec{\Psi}_0 $ is lower triangular, as well as the standard definition of the
operator norm of a matrix in the penultimate step.  We define the {\em Stochastic Non-Isotropic Finite-Time Lyapunov Exponent}
(SNIFTLE, pronounced `sniffle') by
\begin{equation}
\lambda \left( \vec{\xi}_0,t; \vec{\Xi}_0 \right) := \frac{1}{t} \ln \left\| \vec{\nabla} \F{0}{t} \vec{\Psi}_0 \right\| \, ,
\label{eq:sniftle}
\end{equation}
which is {\em probabilistic} since it is associated with an initial condition in $ {\mathcal N} \left( \vec{\xi}_0, \delta^2
\vec{\Psi}_0 \vec{\Psi}_0^T \right) $.   SNIFTLE clearly generalizes the standard FTLE \eqref{eq:ftledef} since it is stochastic, and allows for a non-isotropic Gaussian probability
distribution with non-infinitesimal extent.  Notably, inserting a radially-symmetric Gaussian  distribution with orthogonal equi-distributional axes by making the choice $ \vec{\Xi}_0 = \vec{I}_n $, the $ n $-dimensional
identity matrix (for which $ \vec{\Psi}_0 = \vec{I}_n $ as well) recovers the {\em deterministic} FTLE. 

Utilizing the general covariance \eqref{eq:covariance}, we note that 
\begin{equation}
Q^2(\vec{\xi}_0,  t) = \frac{1}{\delta^2} \left\| \vec{\Lambda}_t \left( \vec{\xi}_0; 0, \delta \right)  \right\|  \, ,
\label{eq:q2}
\end{equation}
symmetric with the observation \eqref{eq:s2}.  Thus, SNIFTLE can also be written as 
\begin{equation}
 \lambda \left( \vec{\xi}_0, \! t; \vec{\Xi}_0 \right) 
= \frac{1}{t}  \ln  \sqrt{Q^2 \left(\vec{\xi}_0,t \right)}  \, ,
\label{eq:qnoftle}
\end{equation}
clarifying its interpretation \eqref{eq:q2def} as a stochastic deviation from the deterministic expectation.

\section{Conclusions}
We have established how the uncertainty covariance matrix $ \vec{\Lambda}_t $ unifies the FTLE framework with that of stochastic sensitivity.  These are respectively related to the uncertainty in the initial condition and the uncertainty in
the evolving dynamics. Through the definition of SNIFTLE \eqref{eq:sniftle}, this means that the FTLE  emerges in terms of the variance of the probability distribution which
satisfies the Fokker--Planck equation \eqref{eq:fokkerplanck}.  We expect these novel understandings to facilitate new
interpretations of old concepts in evolving physical systems with uncertainty.

\begin{acknowledgments}
{\bf Acknowledgements:} SB acknowledges partial support from the Australian Research Council (grant DP200101764).
\end{acknowledgments}

%

\end{document}